\documentclass[12pt]{article}
\usepackage{pdfsync}
\usepackage{amsmath}
\usepackage{amsfonts}
\usepackage{amssymb}
\usepackage{amsthm}
\usepackage{hyperref}
\usepackage{enumerate}
\usepackage{url}
\usepackage{tikz}
\usepackage{pdflscape}
\usepackage{longtable}
\usepackage{booktabs}
\usepackage{colortbl}
\newcommand{\evnrow}{\rowcolor[gray]{0.95}}
\newcommand{\oddrow}{}

\newcommand{\C}{\mathbb C}
\renewcommand{\P}{\mathbb P}
\newcommand{\wP}{w\mathbb P}
\newcommand{\Q}{\mathbb Q}
\newcommand{\Z}{\mathbb Z}
\newcommand{\N}{\mathbb N}
\newcommand{\cB}{\mathcal{B}}
\newcommand{\cO}{\mathcal{O}}
\newcommand{\cX}{\mathcal{X}}
\newcommand{\cY}{\mathcal{Y}}
\newcommand{\cT}{\mathcal{T}}

\newcommand{\Ygen}{Y_{\mathrm{gen}}}
\newcommand{\Pnum}{P_{\mathrm{num}}}
\newcommand{\broken}{\dasharrow}
\DeclareMathOperator{\wt}{wt}

\DeclareMathOperator{\Pf}{Pf}
\DeclareMathOperator{\Proj}{Proj}

\DeclareMathOperator{\GL}{{GL}}
\DeclareMathOperator{\Hom}{Hom}
\DeclareMathOperator{\Pic}{Pic}
\DeclareMathOperator{\Jer}{Jer}
\DeclareMathOperator{\Tom}{Tom}
\newcommand{\Pini}{P_{\mathrm{ini}}}
\newcommand{\Porb}{P_{\mathrm{orb}}}

\newcommand{\PxP}{\P^2\times\P^2}
\newcommand{\Si}{\Sigma}

\newcommand{\grdb}{{\sc{Grdb}}}

\newtheorem{thm}{Theorem}

\newtheorem{lem}[thm]{Lemma}
\newtheorem{prop}[thm]{Proposition}

\theoremstyle{definition}

\theoremstyle{remark}

\theoremstyle{remark}
\newtheorem{eg}[thm]{Example}

\numberwithin{equation}{section}
\numberwithin{thm}{section}

\newcommand{\QED}{\ifhmode\unskip\nobreak\fi\quad\ensuremath{\mathrm{QED}}}
\newenvironment{pf}{\paragraph{Proof}}{\QED\medskip}

\title{Fano 3-folds in $\P^2\times\P^2$ format, Tom and Jerry}
\author{Gavin Brown, Alexander M. Kasprzyk, Muhammad Imran Qureshi}
\date{}

\begin{document}

\maketitle

\begin{abstract}
We study $\Q$-factorial terminal Fano 3-folds whose equations
are modelled on those of the Segre embedding of $\PxP$.
These lie in codimension~4 in their total anticanonical embedding
and have Picard rank 2. They fit into the
current state of classification in three different ways.
Some families arise as unprojections of degenerations
of complete intersections, where the generic unprojection
is a known prime Fano 3-fold in codimension~3; these
are new, and an analysis of their Gorenstein projections
reveals yet other new families.
Others represent the ``second Tom'' unprojection families
already known in codimension~4, and we show that every such
family contains one of our models.
Yet others have no easy Gorenstein projection analysis at all,
so prove the existence of Fano components on their Hilbert scheme.
\end{abstract}

\section{Introduction}

\subsection{Fano 3-folds, Gorenstein rings and $\PxP$}

A {\em Fano 3-fold} is a complex projective variety $X$ of dimension~3
with $\Q$-factorial terminal singularities and $-K_X$ ample.
We construct several new Fano 3-folds, and others
which explain known phenomena.
The anticanonical ring $R(X) = \oplus_{m\in\N} H^0(X,-mK_X)$
of a Fano 3-fold $X$ is Gorenstein, and
provides an embedding $X\subset \wP$ in weighted projective space (wps)
that we exploit here, focusing on the case $X\subset \wP^7$ of codimension~4.

According to folklore, when seeking Gorenstein rings in codimension~4
one should look to $\P^2\times\P^2$ and $\P^1\times\P^1\times\P^1$.
Each embeds by the Segre embedding
as a projectively normal variety in codimension~4
with Gorenstein coordinate ring (by \cite[\S5]{GW}
since their hyperplane sections are subcanonical).
We consider $W = \PxP$, expressed as
\begin{equation}\label{eq!p2xp2}
W \buildrel{\cong}\over{\longrightarrow} V \colon
\left( \bigwedge^2
\begin{pmatrix} x_1 & x_2 & x_3 \\ x_4 & x_5 & x_6 \\ x_7 & x_8 & x_9 \end{pmatrix}
= 0 \right)
\subset\P^8,
\end{equation}
or, in words, as the locus where a generic $3\times 3$ matrix
of forms drops rank. As part of a more general theory of weighted
homogeneous varieties, the case of $\PxP$ was worked out
by Szendr{\H o}i \cite{sz05}, which
was the inspiration for our study here.

The number of deformation families of Fano 3-folds is finite \cite{K92,kmm},
and the Graded Ring Database (\grdb) \cite{grdb,k3db} has a list of rational functions $P(t)$
that includes all Hilbert series $P_{X}(t) = \sum_{m\in\N} h^0(-mK_X) t^m$
of Fano 3-folds with $\Pic(X)=\Z\cdot ({-}K_X)$.
(In fact, we do not know of
any Fano 3-fold whose Hilbert series is not on that list, even without
this additional condition.)
An attempt at an explicit classification, outlined in \cite{ABR02},
aims to describe all deformation families of Fano 3-folds for
each such Hilbert series.
All families whose general member lies in codimension~$\le2$
are known \cite{ccc}, and almost certainly those in codimension~3 are too \cite{ABR02,grdb}.
An analysis of (Gorenstein) projections \cite{TJ,Pap,RT}
provides much of the classification in codimension~4, but it is not
complete, and codimension~4 remains at the cutting edge.

We use the methods of \cite{TJ} freely, although we work through
an example in detail in \S\ref{s!26989} and explain any novelties as they arise.

\subsection{The aims of this paper}
We describe families of Fano 3-folds $X\subset\wP^7$
whose equations are a specialisation of the format \eqref{eq!p2xp2};
that is, they are regular pullbacks, as in \S\ref{s!V}.
It is usually hard to describe the equations of varieties in
codimension~4---see papers from Kustin and Miller \cite{KM} to
Reid \cite{reid4}---but if we decree the format in advance, then
the equations come almost for free, and the question becomes how to put
a grading on them to give Fano 3-folds.
Our results come in three broad flavours, which we explain in \S\S\ref{s!codim3}--\ref{s!noproj}
and summarise here.

\paragraph{\S\ref{s!codim3} Unprojecting Pfaffian degenerations:}
We find new varieties in $\PxP$ format
that have the same Hilbert series as known Fano 3-folds
but lie in different deformation families.
From another point of view, we understand this as the unprojection
analysis of degenerations of complete intersections, and
this treatment provides yet more families not exhibited by \cite{TJ}.
(The key point is that the unprojection divisor $D\subset Y$ does
not persist throughout the degeneration $Y\leadsto Y_0$, and so the resulting
unprojection is not a degeneration in a known family.)

For example, No.~1.4 in Takagi's analysis \cite{takagi02} exhibits
a single family of Fano 3-folds with Hilbert series
\[
P_{26989}(t) = \frac{1 - 3t^2 - 4t^3 + 12t^4  - 4t^5 - 3t^6 + t^8}{(1-t)^7(1-t^2)}
  = 1 + 7t + 26t^2 + 66t^3 + \cdots;
\]
this is number 26989 in the \grdb.
Our $\PxP$ analysis finds another family with $\rho_X=2$, and a
subsequent degeneration--unprojection
analysis of the situation finds a third family.
\begin{thm}
\label{th!712}
There are three deformation families of Fano 3-folds~$X$
with Hilbert series $P_X = P_{26989}$.
Their respective general members $X\subset \P(1^7,2)$
all lie in codimension~4 with degree $-K_X^3 = 17/2$ and
a single orbifold singularity $\frac{1}2(1,1,1)$, and
with invariants:
\[
\begin{array}{l|ccc|lc}
& \rho_X & e(X) & h^{2,1}(X) & \textrm{Construction}&N\\
\hline
\textrm{Family 1 (\cite[1.4]{takagi02})} & 1 & -14 & 9 & \S\ref{sec!special0}: c.i.\ unprojection &6 \\
\textrm{Family 2} & 2 & -16 & 11 & \S\ref{sec!special1}: \Tom_3&5 \\
\textrm{Family 3} & 2^* & -12 & 9^* & \S\ref{sec!special2}: \Jer_{1,3} &7
\end{array}
\]
(The superscript ${}^*$ in Family~3 indicates a computer algebra calculation.)
\end{thm}
We prove this particular result in \S\ref{s!26989};
the last two columns of the table refer to the unprojection calculation
($N$ is the number of nodes, as described in \S\ref{s!26989}), which
is explained in the indicated sections.
The Euler characteristic $e(X)$ is calculated during the unprojection
following \cite[\S7]{TJ} and the other invariants follow.
We do not know whether there are any other deformation families
realising the same Hilbert series $P_X=P_{26989}(t)$.

We calculate the Hodge number $h^{2,1}(X)$ in Family~3 using
Ilten's computer package \cite{ilten} for the computer algebra system
Macaulay2 \cite{M2} following
\cite{DFF}: denoting the affine cone over $X$ by $A_X$,
Theorem~2.5 of~\cite{DFF} gives
\[
H^{2,1}(X) \cong \left(T_{A_X}^1\right)_{-1},
\]
and this is exactly what \cite{ilten} calculates (compare \cite[\S4.1.3]{BF}).

In this case, all three families lie in codimension~4. It is more
common that the known family lies in codimension~3
and we find new families in codimension~4.
Thus the corresponding Hilbert scheme contains different components
whose general members are Fano 3-folds in different codimensions,
a phenomenon we had not seen before.

Further analysis of degenerations finds yet more new Fano 3-folds
even where there is no $\PxP$ model;
the following result is proved in \S\ref{s!degen}.

\begin{thm}
\label{th!548}
There are two deformation families of Fano 3-folds~$X$
with Hilbert series $P_X = P_{548}$.
Their respective general members $X$ have degree $-K_X^3 = 1/15$,
and are distinguished by their embedding in wps and Euler characteristics
as follows:
\[
\begin{array}{l|cc|cc}
& X\subset\wP & e(X) & \# nodes\\
\hline
\textrm{Family 1} & X\subset\P(1,3,4,5,6,7,10) & -42 & 8 \\
\textrm{Family 2} & X\subset\P(1,3,4,5,6,7,9,10) & -40 &9 
\end{array}
\]
\end{thm}
\noindent
In this case there is no $\PxP$ model: such a model would come
from a specialised Tom unprojection, but the Tom and Jerry
analysis outlined in \S\ref{s!degen} rules this out.

\paragraph{\S\ref{s!tom} Second Tom:}
The Big Table \cite{BKRbigtable} lists all (general) Fano 3-folds
in codimension~4 that have a Type~I projection. Such projections
can be of Tom type or Jerry type (see \cite[2.3]{TJ}). The result of
that paper is that every Fano 3-fold admitting a Type~I projection
has at least one Tom family and one Jerry family. However in some cases
there is a second Tom or second Jerry (or both).
Two of these cases were already known to Szendr{\H o}i \cite{sz05},
even before the Tom and Jerry analysis was developed.

Euler characteristic is of course constant in families, but
whenever there is a second Tom, the Euler characteristics of
members of the two Tom families differ by~2.
Theorem~\ref{th!secondtom} below says that in this case the Tom family with smaller
Euler characteristic always contains special members  in $\PxP$ format.

\paragraph{\S\ref{s!noproj} No Type I projection:}
Finally, we find some Fano 3-folds that are harder to describe,
including some that currently have no construction by Gorenstein unprojection.
Such Fano 3-folds were expected to exist, but this is the first construction
of them in the literature we are aware of. It may be the case that there are
other families of such Fano 3-folds having Picard rank~1, but our
methods here cannot answer that question.

\subsection{Summary of results}

Our approach starts with a systematic enumeration of all possible
$\PxP$ formats that could realise the Hilbert series of a Fano 3-fold
after appropriate specialisation. 
In \S\ref{s!formats},
following~\cite{BKZ,QJSC}, we find 53 varieties in $\PxP$ format
that have the Hilbert series of a Fano 3-fold.
We summarise the fate of each of these 53 cases in Table~\ref{tab!results};
the final column summarises our results, as we describe below, and the
rest of the paper explains the calculations that provide the proof.

\begin{longtable}{>{\hspace{0.5em}}llcccccl<{\hspace{0.5em}}}
\caption{53 Fano 3-fold Hilbert series in $\PxP$ format. (Number of nodes is given as a superscript to $\Tom$/$\Jer$.)} \label{tab!results}\\
\toprule
\multicolumn{1}{c}{$k$}&\multicolumn{1}{c}{$a$}&\multicolumn{1}{c}{$b$}&\multicolumn{1}{c}{\grdb}&\multicolumn{1}{c}{$c$}&\multicolumn{1}{c}{T/J}&\multicolumn{1}{c}{$\wP$ in \grdb}&\multicolumn{1}{c}{codim 4 models in this paper}\\
\cmidrule(lr){1-3}\cmidrule(lr){4-4}\cmidrule(lr){5-7}\cmidrule(lr){8-8}
\endfirsthead
\multicolumn{8}{l}{\vspace{-0.25em}\scriptsize\emph{\tablename\ \thetable{} continued from previous page}}\\
\midrule
\endhead
\multicolumn{8}{r}{\scriptsize\emph{Continued on next page}}\\
\endfoot
\bottomrule
\endlastfoot
\evnrow 4 & 000 & 112 & 26989 & 4 &  & $\P(1^7,2)$  &  $\Tom^5,$ $\Jer^7$ in $\P(1^7,2)$\\
\oddrow 5 & 000 & 122 & 20652 & 4  & TTJ\ \ & $\P(1^5,2^3)$ & second Tom \\
\evnrow  5 & 001 & 112 & 20543 & 3 & n/a & $\P(1^5,2^2)$ &$\Tom^7$, $\Jer^9$ in $\P(1^5,2^3)$ \\
5 & 001 & 112& 24078 & 4 & TTJ\ \ & $\P(1^6,2,3)$ & second Tom \\
\evnrow 6 & 000 & 222 & 12960 & 4 &\ TJ& $\P(1^3,2^5)$ & subfamily of Tom \\
6&001&122&16339&4&TTJJ&$\P(1^4,2^3,3)$&second Tom\\
\evnrow 7&001&123&11436&3&n/a&$\P(1^3,2^3,3)$&$\Tom^{13}$ in $\P(1^3,2^3,3^2)$\\
7&001&123&16228&4&TTJJ&$\P(1^4,2^2,3,4)$&second Tom\\
\evnrow 7&011&122&11455&4&TTJJ&$\P(1^3,2^3,3^2)$&second Tom\\
8&001&223&11157&5&n/a&$\P(1^3,2^2,3^2,4^2)$&bad $1/4$ point\\
\evnrow 8&001&223&6878&4&TTJJ&$\P(1^2,2^3,3^3)$&second Tom\\
8&011&123&11125&4&TTJJ&$\P(1^3,2^2,3^2,4)$&second Tom\\
\evnrow 9&001&233&5970&4&TTJJ&$\P(1^2,2^2,3^3,4)$&second Tom\\
9&012&123&11106&4&TTJJ&$\P(1^3,2^2,3,4,5)$&second Tom\\
\evnrow 9&012&123&11021&4&TTJJ&$\P(1^3,2,3^2,4^2)$&second Tom\\
9&012&123&5962&3&n/a&$\P(1^2,2^2,3^3)$&$\Tom^{11}$, $\Jer^{13}$ in $\P(1^2,2^2,3^3,4)$\\
\evnrow 9&012&123&6860&4&TTJ\ \ &$\P(1^2,2^3,3^2,5)$&second Tom\\
10&001&234&5870&4&TTJJ&$\P(1^2,2^2,3^2,4,5)$&second Tom\\
\evnrow 10&011&233&5530&4&TTJJ&$\P(1^2,2,3^3,4^2)$&second Tom\\
10&012&124&10984&3&n/a&$\P(1^3,2,3,4,5)$&bad 1/4 point\\
\evnrow 10&012&124&5858&3&n/a&$\P(1^2,2^2,3^2,5)$&$\Tom^{13}$, $\Jer^{14}$ in $\P(1^2,2^2,3^2,4,5)$\\
11&011&234&5306&4&TTJJ&$\P(1^2,2,3^2,4^2,5)$&second Tom\\
\evnrow 11&012&134&5302&3&n/a&$\P(1^2,2,3^2,4^2)$&$\Tom^{16}$ in $\P(1^2,2,3^2,4^2,5)$\\
11&012&134&5844&3&n/a&$\P(1^2,2^2,3,4,5)$&bad $1/6$ point and no $1/5$\\
\evnrow 11&012&134&10985&4&TTJJ&$\P(1^3,2,3,4,5,6)$&second Tom\\
12&012&234&1766&4&no I&$\P(1,2,3^3,4^2,5)$&quasismooth $\PxP$ model\\
\evnrow 12&012&234&5215&4&TTJJ&$\P(1^2,2,3,4^2,5^2)$&second Tom\\
12&012&234&2427&4&TTJJ&$\P(1,2^2,3^2,4,5^2)$&second Tom\\
\evnrow 12&012&234&5268&4&TTJJ&$\P(1^2,2,3^2,4,5,6)$&second Tom\\
13&001&345&1413&4&TTJJ&$\P(1,2,3^2,4^2,5^2)$&second Tom\\
\evnrow 13&012&235&5177&4&TJ&$\P(1^2,2,3,4,5^2,6)$&bad 1/5 point\\
13&012&235&2422&4&TTJJ&$\P(1,2^2,3^2,4,5,7)$&second Tom\\
\evnrow 14&011&345&5002&4&TTJJ&$\P(1^2,3,4^2,5^2,6)$&second Tom\\
14&012&245&5163&4&TTJJ&$\P(1^2,2,3,4,5,6,7)$&second Tom\\
\evnrow 14&012&245&1410&4&\ \ TJJ&$\P(1,2,3^2,4^2,5,7)$&bad 1/4 point\\
14&013&235&4999&3&n/a&$\P(1^2,3,4^2,5^2,6)$&bad 1/4 point\\
\evnrow 14&013&235&1396&3&n/a&$\P(1,2,3^2,4,5^2)$&$\Tom^9$, $\Jer^{11}$ in $\P(1,2,3^2,4,5^2,6)$\\
15&012&345&878&4&no I&$\P(1,3^2,4^2,5^2,6)$&quasismooth $\PxP$ model\\
\evnrow 15&012&345&4949&4&TTJJ&$\P(1^2,3,4,5^2,6^2)$&second Tom\\
15&012&345&1253&4&TTJ\ \ &$\P(1,2,3,4^2,5^2,7)$&second Tom\\
\evnrow 15&012&345&1218&4&TTJJ&$\P(1,2,3,4,5^3,6)$&second Tom\\
15&012&345&4989&4&TTJJ&$\P(1^2,3,4^2,5,6,7)$&second Tom\\
\evnrow 16&012&346&1186&4&TJJ\ \ &$\P(1,2,3,4,5^2,6,7)$&bad 1/5 point\\
17&012&356&648&4&no I&$\P(1,3,4^2,5^2,6,7)$&bad 1/5 point\\
\evnrow 17&012&356&4915&4&TTJJ&$\P(1^2,3,4,5,6,7,8)$&second Tom\\
18&012&456&577&4&no I&$\P(1,3,4,5^2,6^2,7)$&quasismooth but not terminal\\
\evnrow 18&012&456&645&4&TJ&$\P(1,3,4^2,5,6,7^2)$&bad $1/4$ point\\
18&012&456&4860&4&TTJJ&$\P(1^2,4,5,6^2,7^2)$&second Tom\\
\evnrow 19&012&457&570&4&\ \ TJJ&$\P(1,3,4,5^2,6,7,8)$&bad 1/5 point\\
20&012&467&4839&4&TTJJ&$\P(1^2,4,5,6,7,8,9)$&second Tom\\
\evnrow 22&012&568&1091&4&\ \ TJJ&$\P(1,2,5,6,7^2,8,9)$&bad 1/7 point\\
22&012&568&393&4&TJ&$\P(1,4,5^2,6,7,8,9)$&bad 1/4, 1/5 points\\
\evnrow 23&012&578&360&4&no I&$\P(1,4,5,6,7^2,8,9)$&bad 1/7 point\\
\end{longtable}

The columns of Table~\ref{tab!results} are as follows.
Column $k$ is an adjunction index, described in \S\ref{s!enum},
and columns $a$ and $b$ refer to the vectors in \S\ref{s!V}
that determine the weights on the weighted $\PxP$.
Column \grdb\ lists the number of the Hilbert series in \cite{grdb},
column $c$ indicates the codimension of the usual model suggested there,
and $\wP$ its ambient space. Column T/J shows the number of distinct
Tom and Jerry components according to \cite{TJ}.
For example, TTJ indicates there are 2 Tom unprojections and
1 Jerry unprojection in the Big Table \cite{BKRbigtable}.
We write `no I' when the Hilbert series does not admit a numerical Type~I projection,
and so the Tom and Jerry analysis does not apply, and `n/a' if the usual model
is in codimension~3 rather than~4.

The final column describes the results of this paper;
it is an abbreviation of more detailed results.
For example, Theorem~\ref{th!712} expands out the first line of the table, $k=4$,
and other lines of the table that are not indicated as failing have analogous theorems that
the final column summarises.
If the $\PxP$ model fails to realise a Fano 3-fold
at all, it is usually because the general member does not
have terminal singularities; we say, for example,
`bad 1/4 point' if the format forces a non-quasismooth, non-terminal
index~4 point onto the variety.

When the \grdb\ model is in codimension~3, we list which Tom and
Jerry unprojections of a degeneration work to give alternative varieties
in codimension~4, indicating the number
of nodes as a superscript and the codimension~4 ambient space.
(We don't say which Tom or Jerry since
that depends on a choice of rows and columns.)
In each case the Tom unprojection gives the $\PxP$ model determined
by the parameters $a$ and $b$.
The usual codimension~3 model
arises by Type~I unprojection with number of nodes being one more
that that of the $\PxP$ Tom model.

When the \grdb\ model is in codimension~4 with 2 Tom unprojections,
the $\PxP$ always works to give the second of the Tom families.
The further Tom and Jerry analysis of the unprojection is carried
out in \cite{TJ} and we do not repeat the result here.
When the \grdb\ model is in codimension~4 with only a single Tom unprojection,
the model usually fails. The exception is family 12960, which does work
as a $\PxP$ model.
There is also a case of a Hilbert series, number 11157, where the \grdb\ offers a prediction
of a variety in codimension~5, but this fails as a $\PxP$ model.

In \S\ref{s!enum}, we outline a computer search that provides the
$a,b$ parameters of Table~\ref{tab!results} which are the starting point
of the analysis here. In \S\ref{s!sz05}, we summarise the results of
\cite{sz05} that provide the most general form of the
Hilbert series of a variety in $\PxP$ format; that paper also discovered
cases 11106 and 11021 of Table~\ref{tab!results} that inspired our approach here.
First we introduce the key varieties of the $\PxP$ format in~\S\ref{s!V}.

\section{The key varieties and weighted $\P^2\times\P^2$ formats}
\label{s!V}

The affine cone $C(\PxP)$ on $\PxP$ is defined by the equations~\eqref{eq!p2xp2}
on $\C^9$.
It admits a \hbox{6-dimensional} family of $\C^*$ actions, or equivalently
6 degrees of freedom in assigning positive integer gradings to
its (affine) coordinate ring. We express this as follows.

Let $a=(a_1,a_2,a_3)$ and $b=(b_1,b_2,b_3)$ be two vectors of integers
that satisfy $a_1\le a_2\le a_3$, and similarly for the $b_i$, and that $a_1+b_1\ge1$.
We define a {\em weighted $\P^2\times\P^2$} as
\begin{equation}\label{eq!Veqns}
V = V(a,b) = \left(
\bigwedge^2
\begin{pmatrix} x_0 & x_1 & x_2 \\ x_3 & x_4 & x_5 \\ x_6 & x_7 & x_8 \end{pmatrix}
= 0 \right) \subset \P^8(a_1+b_1,\dots,a_3+b_3),
\end{equation}
where the variables have weights
\begin{equation}\label{eq!wtmx}
\wt
\begin{pmatrix} x_0 & x_1 & x_2 \\ x_3 & x_4 & x_5 \\ x_6 & x_7 & x_8 \end{pmatrix}
=
\begin{pmatrix} a_1+b_1 & a_1+b_2 & a_1+b_3 \\ a_2+b_1 & a_2+b_2 & a_2+b_3 \\ a_3+b_1 & a_3+b_2 & a_3+b_3  \end{pmatrix}
=: a^T + b.
\end{equation}
Thus $V(a,b) = C(\PxP) /\!\!/ \C^*$, where the $\C^*$ action is determined
by the grading.
We treat $V(a,b)$ as a key variety for each different pair $a,b$.
(Note that the entries of $a$ and $b$ may also all lie in $\frac{1}2+\Z$,
without any change to our treatment here.)

\begin{prop}
$V(a,b)$ is a 4-dimensional, $\Q$-factorial projective toric variety of
Picard rank $\rho_V = 2$.
\end{prop}

\begin{pf}
First we describe a toric variety $W(a,b)$ by its Cox ring.
The input data is the weight matrix \eqref{eq!wtmx}, which is weakly increasing along rows
and down columns.
The key is to understand the freedom one has to choose alternative vectors
$a^{(i)}$, $b^{(i)}$, for $i = 1,2$, to give the same matrix.
For example, if we choose $a^{(1)}_1=0$, then $b^{(1)}$ is determined by the top row, and
then $a^{(1)}_2$ and $a^{(1)}_3$ are determined by the first column.
Alternatively, choosing $b^{(2)}_1=0$ determines different vectors $a^{(2)}$ and $b^{(2)}$.
Concatenating the $a$ and $b$ vectors to give $v^{(i)}=(a^{(i)}_1,\dots,b^{(i)}_3)\in\Q^6$
determines a 2-dimensional $\Q$-subspace $U=U_{a,b}\subset\Q^6$
together with a chosen integral basis $\left<v^{(1)}, v^{(2)}\right>$.

We define $W(a,b)$ as a quotient of $\C^6$ by $\C^*\times\C^*$ as follows.
In terms of Cox coordinates, it is determined by the polynomial ring $R$
in variables $u_1,u_2,u_3$, $v_1,v_2,v_3$,
bi-graded by the columns of the matrix (giving the two $\C^*$ actions)
\begin{equation}
\label{e!wtmx}
\begin{pmatrix}
a^{(1)}_1 & a^{(1)}_2 & a^{(1)}_3 & b^{(1)}_1 & b^{(1)}_2 & b^{(1)}_3 \\
a^{(2)}_1 & a^{(2)}_2 & a^{(2)}_3 & b^{(2)}_1 & b^{(2)}_2 & b^{(2)}_3
\end{pmatrix}.
\end{equation}
The irrelevant ideal is $B(a,b) = \left< u_1,u_2.u_3 \right> \cap \left< v_1,v_2.v_3 \right>$, and
\[
W(a,b) = \left(\C^6 \setminus V(B(a,b))\right) / \ \C^*\times\C^*.
\]
If $W(a,b)$ is well formed, then it is a toric variety determined by a fan
(the image of all non-irrelevant cones of the fan of $\C^6$ under projection
to a complement of~$U$).
The bilinear map
\begin{eqnarray}\label{eq!phi}
\begin{array}{rrcl}
\Phi_{a,b} \colon & W(a,b) & \longrightarrow & \P(a_1+b_1,a_1+b_2,\dots,a_3+b_3) \\
 & (u_1,\dots,v_3) & \mapsto & (u_1v_1,u_1v_2,\dots,u_3v_3)
\end{array}
\end{eqnarray}
is an isomorphism onto its image $V(a,b)$,
and the conclusions of the proposition all follow at once.
($\Q$-factoriality holds since the Cox coordinates correspond to the 1-skeleton
of the fan, and so any maximal cone with at least 5 rays must contain
all $u_i$ or all $v_j$, contradicting the choice of irrelevant ideal.)

If $W(a,b)$ is not well formed, then, just as for wps
(see Iano-Fletcher \cite[6.9--20]{Fletcher}), there is a different weight matrix
that is well formed and determines a toric variety $W'$ isomorphic to $W(a,b)$.
The proposition follows using $W'$.
\end{pf}

We review the well forming process used in the proof. It has two parts.
If some integer $n$ divides every entry of the
weight matrix \eqref{e!wtmx} then we may divide through by that.
The subspace $U\subset\Q^6$ is unchange by this; the grading on $R$
is simply scaled by~$n$.
Otherwise some integer $n$ divides all columns except one.
In that case, the corresponding Cox coordinate $u$ appears only
as $u^n$ in the coordinate rings of standard affine patches.
We may truncate $R$ by replacing the generator $u$ by $u^n$; this
does not change the coordinate rings of the affine patches, and so
the scheme it defines is isomorphic to the original.
This multiplies the $u$ column of \eqref{e!wtmx} by $n$, changing
the subspace $U$, and then we may divide the whole matrix by $n$
as before. 

Having said that, in practice we will work with non-well-formed quotients
if they arise, since they still admit regular pullbacks that are well formed,
and the grading on the target wps is something we fix in advance.
More importantly for us here is that well forming step $u\leadsto u^n$
destroys the $\PxP$ structure, so we avoid it.

\begin{eg}
Consider $V(a,b) \subset \P(2^6,3^3)$ for $a=(1,1,1)$, $b=(1,1,2)$.
Selecting $a^{(i)}$ and $b^{(i)}$ as above gives bi-grading matrix
\[
\left(
\begin{array}{ccc|ccc}
0 & 0 & 0 & 2 & 2 & 3 \\
2 & 2 & 2 & 0 & 0 & 1
\end{array}
\right)
\]
on variables $u_1,u_2,u_3$, $v_1,v_2,v_3$.
(We use the vertical line in the bi-grading matrix to indicate the
irrelevant ideal $B(a,b)$.)
The map $\Phi$ of \eqref{eq!phi} is then
\begin{eqnarray*}
W(a,b) & \longrightarrow & \P(2,2,3,2,2,3,2,2,3) = \P(2^6,3^3) \\
(u_1,\dots,v_3) & \mapsto & (u_1v_1, u_1v_2,\dots,u_3v_3),
\end{eqnarray*}
since the monomials having gradings
$\left(\begin{smallmatrix}2\\2\end{smallmatrix}\right)$ and
$\left(\begin{smallmatrix}3\\3\end{smallmatrix}\right)$,
as necessary.
The image $V(a,b)$ is defined by \eqref{eq!Veqns},
and we often write the target weights of $\Phi$ in matching array:
\[
\begin{pmatrix} 2&2&3 \\2&2&3 \\2&2&3
\end{pmatrix}.
\]

In this case $V(a,b)$ is not well formed: the locus $V(a,b) \cap \P(2^6)$
has dimension~3 (by Hilbert--Burch), so has codimension~1 in $V(a,b)$
but nontrivial stabiliser $\Z/2$ in the wps.
Well forming the gradings using $v_3^2$, as above, gives a new bi-grading
\[
\left(
\begin{array}{ccc|ccc}
0 & 0 & 0 & 1 & 1 & 3 \\
1 & 1 & 1 & 0 & 0 & 1
\end{array}
\right).
\]
That process is well established,
but has a problem: for this presentation $W'$ of $W$, the Segre map is not bi-linear:
$u_1v_1$ has bidegree
$\left(\begin{smallmatrix} 1\\1 \end{smallmatrix}\right)$,
but $u_1v_3$ has an independent bidegree
$\left(\begin{smallmatrix} 3\\2 \end{smallmatrix}\right)$.
We could use $u_1^2v_3$ instead, which has proportional bidegree 
$\left(\begin{smallmatrix} 3\\3 \end{smallmatrix}\right)$.
Taking $V'=\Proj R$, where $R$ is the graded ring of forms
of degrees
$\left(\begin{smallmatrix} m\\m \end{smallmatrix}\right)$
for $m\ge0$, gives $W'\rightarrow V'\subset\P(1^6,3^6)$,
which is now well formed, but we have lost the codimension~4
property of $V$ we want to exploit. In a case like this,
we work directly with the non-well-formed $W(a,b)$
and its non-well-formed image $V\subset\P(2^6,3^3)$.
\end{eg}


We use the varieties $V(a,b)$ as key varieties to produce
new varieties from by regular pullback; see \cite[\S1.5]{fun}
or \cite[\S2]{BKZ}.
In practical terms, that means writing equations in the
form of \eqref{eq!p2xp2} inside a wps $w\P^7$ where
the $x_i$ are homogeneous forms of positive degrees,
and the resulting loci $X\subset w\P^7$ are the
Fano 3-folds we seek.

Alternatively, we may treat $X$ as a complete intersection
in a projective cone over $V(a,b)$, as in~\ref{sec!special1} below,
where the additional cone vertex variables may have any positive degrees;
this point of view is taken by Corti--Reid and Szendr{\H o}i in \cite{CR,sz05,QS,qureshi}.
It follows from this description that the Picard rank of $X$ is~2.

\section{Unprojection and the proof of Theorem~\ref{th!712}}
\label{s!26989}

The Hilbert series number 26989 in the Graded Ring Database (\grdb) \cite{grdb} is
\[
P = \frac{1 - 3t^2 - 4t^3 + 12t^4  - 4t^5 - 3t^6 + t^8}{(1-t)^7(1-t^2)}.
\]
In \S\ref{sec!special0} we describe the known family of Fano 3-folds
$X^{(1)}\subset\P(1^7,2)$  that realise this Hilbert series, $P_{X^{(1)}} = P$.
These 3-folds are not smooth: 
the general member of the family has a single $\frac{1}2(1,1,1)$ quotient singularity.
We exhibit a different family in \S\ref{sec!special1}
with the same Hilbert series in $\PxP$ format, and
the subsequent ``Tom and Jerry'' analysis yields a third distinct family
in \S\ref{sec!special2}.

Recall (from \cite[\S4]{TJ}, for example)
that if $X\broken Y\supset D$ is a Gorenstein unprojection and $Y$
is quasismooth away from $N$ nodes, all of which lie on $D$, then
\begin{equation}
\label{eq!nodes}
e(X) = e(Y) + 2N - 2.
\end{equation}

\subsection{The classical $7\times 12$ family}
\label{sec!special0}

A general member of the first family
can be constructed as the unprojection of a coordinate
$D = \P^2$ inside a c.i.\ $Y_{2,2,2}\subset\P^6$ (see, for example, Papadakis \cite{Pwith}).
In general, $Y$ has 6 nodes that lie on $D$: in coordinates $x,y,z$,
$u,v,w,t$ of $\P^6$, setting $D = (u=v=w=t=0)$, the general $Y$
has equations defined by
\[
\begin{pmatrix}
A_{1,1} & \cdots & A_{1,4} \\ 
A_{2,1} & \cdots & A_{2,4} \\ 
A_{3,1} & \cdots & A_{3,4} 
\end{pmatrix}
\begin{pmatrix}
u\\v\\w\\t
\end{pmatrix}
=
\begin{pmatrix} 0 \\ 0 \\ 0
\end{pmatrix},
\]
for general linear forms $A_{i,j}$; singularities occur when the
$3\times 4$ matrix drops rank, which is calculated by evaluating
the numerator of the Hilbert series of that locus at 1:
\[
P_{\mathrm{sings}} =
\frac{1 - 4t^3 + 3t^4}{(1-t)^3}
=
\frac{1 + 2t + 3t^2}{1-t},
\quad
\text{so there are $1+2+3=6$ nodes}.
\]

The coordinate ring of $X$ has a $7\times 12$ free resolution.
If $\Ygen$ is a nonsingular small deformation of $Y$, then
$e(\Ygen) = -24$ (by the usual Chern class calculation, since $\Ygen$
is a smooth $2,2,2$ complete intersection) so, by \eqref{eq!nodes},
\[
e(X) = -24 + 12-2 = -14.
\]
This family is described by Takagi \cite{takagi02};
it is No.\ 1.4 in the tables there of Fano 3-folds of Picard rank~1.

\subsection{A $\PxP$ family with Tom projection}
\label{sec!special1}

Consider the $\PxP$ key variety
$V_{a,b} \subset \P(1^6,2^3)$,
where $a = (\frac{1}2,\frac{1}2,\frac{1}2)$ and $b=(\frac{1}2,\frac{1}2,\frac{3}2)$.
We define a quasismooth variety $X^{(2)} \subset \P(1^7,2)$ in codimension~4
as a regular pullback.

In explicit terms, in coordinates $x,y,z,t,u,v,w,s$ on $\P(1^7,2)$,
a $3\times 3$ matrix $M$ of forms of degrees
\[
a^T + b =
\begin{pmatrix} 1 & 1 & 2 \\ 1 & 1 & 2 \\ 1 & 1 & 2 \end{pmatrix}
\]
gives a quasismooth $X^{(2)} = (\wedge^2M=0)\subset\P(1^7,2^2)$; for example,
\[
M =
\begin{pmatrix}
    x &          t &      s \\
    y &          u &      x^2 - z^2 + t^2 + v^2 \\
    z &          v &      xt + yu + w^2
\end{pmatrix}
\]
works.
Alternatively, note that $X^{(2)}$ may be viewed as a complete intersection
\[
X^{(2)} = C_1V_{a,b} \cap Q_1 \cap Q_2 \subset \P(1^7,2^3),
\]
where $C_1V_{a,b}\subset\P(1^7,2^3)$ is the projective cone over $V_{a,b}$
on a vertex of degree~1 (by introducing a new variable of degree~1),
and $Q_i$ are general quadrics (which are quasilinear, and so may be used
to eliminate two variables of degree~2).
The general such $X^{(2)}$ is quasismooth (since
in particular the intersection misses the vertex).
Described in these terms, $C_1V_{a,b}$ has Picard rank~2, and
so $\rho_{X^{(2)}} = 2$.

Any such $X^{(2)}$ has a single quotient singularity $\frac{1}2(1,1,1)$,
at the coordinate point $P_s\in X^{(2)}$ as the explicit equations make clear,
since $y, z, u, v$ are implicit functions in a neighbourhood of $P_s\in X^{(2)}$.
The Gorenstein projection from this point $P_s$ has image
$Y = \left(\Pf N  = 0 \right) \subset\P^6$, where
\[
N = 
\begin{pmatrix}
0 & x & y & z \\
   & t & u & v \\
       && x^2 - z^2 + t^2 + v^2 & xt + yu + w^2 \\
        &&& 0
\end{pmatrix}
\]
is an antisymmetric $5\times 5$ matrix, and $\Pf N$ denotes the sequence
of 5 maximal Pfaffians of~$N$. (The nonzero entries of $N$
are those of $M^T$ with the entry $s$ deleted.)

This $Y$ contains the projection divisor $D = (y = z = u = v = 0)$
and has 5 nodes on $D$ (either by direct calculation, or by
the formula of \cite[\S7]{TJ}).
The divisor $D\subset Y$ is in Tom$_3$ configuration: entries
$n_{i,j}$ of the skew $5\times 5$ matrix $N$ defining~$Y$ lie
in the ideal $I_D = (y,z,u,v)$ if both $i\not=3$ and $j\not=3$; that is,
all entries off row~3 and column~3 of $N$ are in~$I_D$.
Thus, in particular, we can reconstruct $X^{(2)}$ from $D\subset Y$ as the Tom$_3$
unprojection.
It follows from Papadakis--Reid \cite[\S2.4]{PR} that $\omega_{X^{(2)}} = \cO_{X^{(2)}}(-1)$
and so $X^{(2)}$ is a Fano 3-fold.

It remains to show that $e(X^{(2)}) = -16$, so that this Fano 3-fold must lie
in a different deformation family from the classical one constructed in~\S\ref{sec!special0}.

The degree of the $(1,2)$ entry $f_{1,2}$ of $N$ is in fact zero while
the degree of $f_{4,5}$ is~2, although each entry is of course the zero polynomial in this case;
we denote this by indicating the degrees of the entries with brackets around
those that are zero in this case:
\begin{equation}
\begin{pmatrix} (0)&1&1&1\\&1&1&1\\&&2&2\\&&&(2) \end{pmatrix}.
\end{equation}
We may deform $Y$ by varying these two entries to $f_{1,2}=\varepsilon$
and $f_{4,5} = \varepsilon f$,
where $\varepsilon\not=0$ and $f$ is a general quadric on $\P^6$
(and, of course, the skew symmetric entries in $f_{2,1}$ and $f_{5,4}$).
Denoting the deformed matrix by $N_\varepsilon$,
and $Y_\varepsilon = (\Pf N_\varepsilon=0)$, we see a small
deformation of $Y$ to a smooth Fano 3-fold $Y_\varepsilon\subset\P^6$
that is a $2,2,2$ complete intersection. (The nonzero constant entries of $N_\varepsilon$
provide two syzygies that eliminate two of the five Pfaffians.)
As in \S\ref{sec!special0}, the smoothing $Y_\varepsilon$ has euler characteristic $-24$,
so by \eqref{eq!nodes} we have that $e_{X^{(2)}} = -24 + 10 - 2 = -16$.

Note that the Pfaffian smoothing $Y_\varepsilon$ of $Y$ destroys the unprojection
divisor $D\subset Y$: for $D$ to lie inside $Y_\varepsilon$ the entries $f_{3,4}$
and $f_{3,5}$ of $N_\varepsilon$ would have to lie in $I_D$
(so $N_\varepsilon$ would be in Jer$_{4,5}$ format with the extra
constraint $f_{4,5}=0$),
but then $Y$ would be
singular along $D$ since 3 of the 5 Pfaffians would lie in $I_D^2$.

\subsection{A third family by Jerry unprojection}
\label{sec!special2}

A Tom and Jerry analysis following \cite{TJ} shows that varieties
$D\subset Y\subset\P^6$ defined by Pfaffians 
as in~\S\ref{sec!special1} by the 
maximal Pfaffians of a syzygy matrix $N$ with weights
\[
\begin{pmatrix} 0 & 1 & 1 & 1 \\ & 1 & 1 & 1 \\ && 2 & 2 \\ &&& 2 \end{pmatrix}.
\]
can also be constructed in Jer$_{1,3}$ format: that is, with
all entries $f_{i,j}$ of $N$ lying in $I_D$ whenever $i$ or $j$ lie in $\{1,3\}$.
The general such $D\subset Y$ has 7 nodes on $D$.
Unprojecting $D\subset Y$ gives a general member $X^{(3)}$ of a third family
with $e(X^{(3)}) = -24 + 2\times 7 - 2 = -12$.

This completes the proof of Theorem~\ref{th!712}.

\section{Unprojecting Pfaffian degenerations}
\label{s!codim3}

\subsection{$\PxP$ models with a codimension 3 Pfaffian component}
\label{s!p2c3}

Each of the Fano Hilbert series
1396, 5302, 5858, 5962, 11436, 20543
is realised by a codimension~3 Pfaffian model, which is the simple default
model presented in the~\grdb. (So too are 4999, 5844 and 10984, but
we do not find new models for these.)
We show that they can also be realised by a $\PxP$ model
in a different deformation family (and sometimes a third model too).
The key point is that a projection of the usual model admits alternative
degenerations in higher codimension that also contain a divisor
that can be unprojected.

For example, consider series number 20543 is
\[
P_{20543}(t) = \frac{1 - 4t^3 + 4t^5 - t^8}{(1-t)^5(1-t^2)^2}.
\]
There is a well-known family that realises this as
$X = \left(\Pf(M) = 0\right) \subset \P(1^5,2^2)$ in codimension~3,
where $M$ has degrees
\[
\begin{pmatrix} 1&1&1&1\\&2&2&2\\&&2&2\\&&&2 \end{pmatrix}.
\]
A typical member of this family has a two $\frac{1}2(1,1,1)$ quotient
singularities, and making the Gorenstein projection from either of them
presents $X$ as a Type~I unprojection of
\[
\P^2 = D \subset Y_{3,3}\subset \P(1^5,2).
\]
In general, $Y$ has $8$ nodes lying on $D$, and it smooths
to a nonsingular Fano 3-fold $\Ygen$ with Euler characteristic $e(\Ygen) = -40$.
Thus a general $X$ has Euler characteristic $e(X) = -40 + 2\times 8 - 2 = -26$.

\paragraph{A quasismooth $\P^2\times\P^2$ family}
We can write another (quasismooth) model $X\subset\P(1^5,2^3)$
in codimension~4 in $\PxP$ format with weights
\[
\begin{pmatrix} 1&1&2\\1&1&2\\2&2&3 \end{pmatrix}.
\]
Projecting from $\frac12(1,1,1)$ has image
$Y = \left(\Pf(M) = 0 \right) \subset\P(1^5,2^2)$
where $M$ has degrees
\begin{equation}
\label{eq!20543}
\begin{pmatrix} (0)&1&1&2\\&1&1&2\\&&2&3\\&&&(3) \end{pmatrix},
\end{equation}
and $Y$ has 7 nodes lying on $D$; in coordinates $x,y,z,t,u,w,v$, we may
take $D=\P^2$ to be $(t=u=v=w=0)$.
By varying the $(1,2)$ entry from zero to a unit, $Y$ has a deformation
to a quasismooth $3,3$ complete intersection $\Ygen$ as before, and so,
$e(X) = e(\Ygen) + 2\times 7 - 2 = -40 +14 - 2 = -28$.
Thus these $\PxP$ models are members of a different deformation
family from the original one.

More is true in this case: the general member of this new deformation family is in $\PxP$ format.
Starting with matrix~\eqref{eq!20543} and $D=\P^2$ as above,
the $(1,2)$ entry of the general Tom$_3$ matrix is necessarily
the zero polynomial. In general, the four entries $(1,4)$, $(1,5)$, $(2,4)$ and $(2,5)$
of the matrix are in the ideal $\left<t,u,v,w\right>$, and for the general member
these four variables are dependent on those entries.
Thus the $(4,5)$ entry can be arranged to be zero
by row-and-column operations.

\paragraph{Another family in codimension 4}
There is a third deformation family in this case.
The codimension~3 format~\eqref{eq!20543}
also admits a Jerry$_{15}$ unprojection with 9 nodes on $D$, giving
$X\subset\P(1^5,2^3)$ in codimension~4 with $e(X) = -24$.

\subsection{Pfaffian degenerations of codimension 2 Fano 3-folds}
\label{s!degen}

The key to the cases in~\S\ref{s!p2c3} that the $\PxP$ model exposes
is the degeneration of a codimension~2 Fano 3-fold.
More generally,
Table~3 of \cite{k3degen} lists 13 cases of Fano 3-fold degenerations
where the generic fibre is a codimension~2 complete intersection
and the special fibre is a codimension~3 Pfaffian.
In each case, the anti-symmetric $5\times 5$ syzygy matrix of the special fibre has
an entry of degree~0, which is the zero polynomial in the degeneration,
but when nonzero serves to eliminate a single variable.
(In fact \cite{k3degen} describes the graded rings of K3 surfaces,
but these extend to Fano 3-folds by the usual extension--deformation
method introducing a new variable of degree~1.)

For example, $Y_{12,13}\subset\P(1,3,4,5,6,7)$ degenerates to
codimension~3
\[
Y^{0} \subset\P(1,3,4,5,6,7,9)
\quad
\text{with syzygy degrees}
\quad
\begin{pmatrix}
0 & 3 & 4 & 7 \\
   & 5 & 6 & 9 \\
     && 9 & 12 \\
         &&& 13
\end{pmatrix}.
\]
Both of these realise Fano Hilbert series number 547, and
the Euler charactistic of a general member is $e(Y) = -56$.

The codimension~2 family has a subfamily whose members
contain a Type~I unprojection divisor,
\[
D = \P(1,3,7) \subset Y = Y_{12,13}\subset \P(1,3,4,5,6,7)
\]
on which $Y$ has 8 nodes;
the unprojection of $D\subset Y$ gives the codimension~3 Pfaffian family
\begin{equation}\label{e!X3}
\text{Hilbert series no.\ 548:}\qquad
X_{12,13,14,15,16}\subset\P(1,3,4,5,6,7,10).
\end{equation}
Imposing the same unprojection divisor $D\subset Y^{0}$ can be
done in two distinct ways, coming from different Tom and Jerry arrangments.
In one way, there are degenerations $Y^{t}_{12,13} \leadsto Y^{0}$
which contain the same $D$ in every fibre $Y_t$. These unproject to
a degeneration of the family \eqref{e!X3} by the following lemma: indeed
unprojection commutes with regular sequences by \cite[Lemma~5.6]{dip1},
and so unprojection commutes with flat deformation, if one fixes the
unprojection divisor; so the lemma is a particular case of  \cite[Lemma~5.6]{dip1}.

\begin{lem}
Let $\P = \P(a_0,\dots,a_s)$ be any wps and
fix $D=\P(a_0,\dots,a_d)\subset\P$, for some $d\le s-2$.
Suppose $Y_t\subset\cY\rightarrow\cT$ is a flat 1-dimensional family of
projectively Gorenstein subschemes of $\P$ over smooth base $0\in\cT$,
each one containing $D$ and with $\dim Y_t = \dim D + 1=d+1$, and with
$\omega_Y=\cO_Y(k_Y)$.
Let $\cX \owns X_t \subset\P(a_0,\dots,a_s,b)$ be
the unprojection of $D\times\cT\subset\cY$,
where $b = k_Y-k_D = a_0+\cdots+a_d-1$.
Then $\cX$ is flat over $\cT$, and for each closed point $t\in\cT$
the fibre $X_t\in\cX$ is the unprojection of $D\subset Y_t$.
\end{lem}

But  the $\Jer_{24}$ unprojection is different: small deformations of
$Y^{0}$ do not contain $D$. Indeed, in this $D\subset Y^0$ model,
$Y^{0}$ has 9 nodes on $D$, which is a numerical obstruction
to any such deformation. This $D\subset Y^{0}$ unprojects
to a codimension~4 Fano 3-fold
\[
X^{0} \subset \P(1,3,4,5,6,7,9,10)
\]
with the same Hilbert series number~548 as \eqref{e!X3}
but lying in a different component: it has Euler characteristic $-56 + 2\times 9 - 2 = -40$.
This proves Theorem~\ref{th!548}.

\section{$\PxP$ and the second Tom}
\label{s!tom}

The Big Table \cite{BKRbigtable}, which contains the results of \cite{TJ},
lists deformation families of Fano 3-folds in codimension~4 that
have a Type~I projection to a Pfaffian 3-fold in codimension~3.
The components are listed according to the Tom or Jerry type
of the projection: the type of projection is invariant for sufficiently
general members of each component.
The result of this section gives an interpretation of the Big Table of \cite{TJ},
but does not describe any new families of Fano 3-folds.

\begin{thm}
\label{th!secondtom}
For every Hilbert series listed in the Big Table \cite{BKRbigtable}
that is realised by two distinct Tom projections,
there is a Fano 3-fold in $\PxP$ format that lies on the family containing
3-folds with the smaller (more negative) Euler characteristic.
\end{thm}

The theorem is proved simply by constructing each case.
There are 29 Hilbert series that have two Tom families.
Using `TTJ'  to indicate a series realised by 2 Tom components and 1 Jerry
component and `TTJJ' to indicate 2 of each, they are:

\begin{longtable}{>{\hspace{0.5em}}cccc<{\hspace{0.5em}}}
\caption{Hilbert series in $\PxP$ format that admit a second Tom unprojection.} \label{tab!resultsraw}\\
\toprule
\multicolumn{1}{c}{\grdb}&\multicolumn{1}{c}{$\PxP$ weights}&\multicolumn{1}{c}{T/J families}&\multicolumn{1}{c}{centre: \#nodes}\\
\cmidrule(lr){1-2}\cmidrule(lr){3-4}
\endfirsthead
\multicolumn{4}{l}{\vspace{-0.25em}\scriptsize\emph{\tablename\ \thetable{} continued from previous page}}\\
\midrule
\endhead
\multicolumn{4}{r}{\scriptsize\emph{Continued on next page}}\\
\endfoot
\bottomrule
\endlastfoot
1253 & $\left(\begin{smallmatrix} 3&4&5\\4&5&6\\5&6&7 \end{smallmatrix}\right)$ & TTJ & $\frac{1}7:6$ \\
1218 & $\left(\begin{smallmatrix} 3&4&5\\4&5&6\\5&6&7 \end{smallmatrix}\right)$ & TTJJ & $\frac{1}5:9$ \\
1413 & $\left(\begin{smallmatrix} 3&4&5\\3&4&5\\4&5&6 \end{smallmatrix}\right)$ & TTJJ & $\frac{1}5:7$ \\
2422 & $\left(\begin{smallmatrix} 2&3&5\\3&4&6\\4&5&7 \end{smallmatrix}\right)$ & TTJJ & $\frac{1}7:5$ \\
2427 & $\left(\begin{smallmatrix} 2&3&4\\3&4&5\\4&5&6 \end{smallmatrix}\right)$ & TTJJ & $\frac{1}5:6$ \\
4839 & $\left(\begin{smallmatrix}
 4 & 6 & 7 \\ 
 5 & 7 & 8 \\
 6 & 8 & 9
\end{smallmatrix}\right)$
& TTJJ & $\frac{1}5:20$; $\frac{1}9 : 13$ \\
4860 & $\left(\begin{smallmatrix} 4&5&6\\5&6&7\\6&7&8 \end{smallmatrix}\right)$ & TTJJ & $\frac{1}7:13$ \\
4915 & $\left(\begin{smallmatrix} 3&5&6\\4&6&7\\5&7&8 \end{smallmatrix}\right)$ & TTJJ & $\frac{1}4:19$; $\frac{1}8:11$ \\
4949 & $\left(\begin{smallmatrix} 3&4&5\\4&5&6\\5&6&7 \end{smallmatrix}\right)$ & TTJJ & $\frac{1}6:11$ \\
4989 & $\left(\begin{smallmatrix} 3&4&5\\4&5&6\\5&6&7 \end{smallmatrix}\right)$ & TTJJ & $\frac{1}4:15$; $\frac{1}7:10$ \\
5002 & $\left(\begin{smallmatrix} 3&4&5\\4&5&6\\4&5&6 \end{smallmatrix}\right)$ & TTJJ & $\frac{1}4:14$; $\frac{1}5:11$; $\frac{1}6:10$ \\
5163 & $\left(\begin{smallmatrix} 2&4&5\\3&5&6\\4&6&7 \end{smallmatrix}\right)$ & TTJJ & $\frac{1}3:19$; $\frac{1}7:9$ \\
5215 & $\left(\begin{smallmatrix} 2&3&4\\3&4&5\\4&5&6 \end{smallmatrix}\right)$ & TTJJ & $\frac{1}5: 9$ \\
5268 & $\left(\begin{smallmatrix} 2&3&4\\3&4&5\\4&5&6 \end{smallmatrix}\right)$ & TTJJ & $\frac{1}3:14$; $\frac{1}5:8$ \\
5306 & $\left(\begin{smallmatrix} 2&3&4\\3&4&5\\3&4&5 \end{smallmatrix}\right)$ & TTJJ & $\frac{1}3:13$; $\frac{1}4:9$; $\frac{1}5:8$ \\
5530 & $\left(\begin{smallmatrix} 2&3&3\\3&4&4\\3&4&4 \end{smallmatrix}\right)$ & TTJJ & $\frac{1}3:11$; $\frac{1}4:8$ \\
5870 & $\left(\begin{smallmatrix} 2&3&4\\2&3&4\\3&4&5 \end{smallmatrix}\right)$ & TTJJ & $\frac{1}3:10$; $\frac{1}5:7$ \\
5970 & $\left(\begin{smallmatrix} 2&3&3\\2&3&3\\3&4&4 \end{smallmatrix}\right)$ & TTJJ & $\frac{1}3:9$; $\frac{1}4: 7$ \\
6860 & $\left(\begin{smallmatrix} 1&2&3\\2&3&4\\3&4&5 \end{smallmatrix}\right)$ & TTJ & $\frac{1}5:4$ \\
6878 & $\left(\begin{smallmatrix} 2&2&3\\2&2&3\\3&3&4 \end{smallmatrix}\right)$ & TTJJ & $\frac{1}3:8$ \\
10985 & $\left(\begin{smallmatrix} 1&3&4\\2&4&5\\3&5&6 \end{smallmatrix}\right)$ & TTJJ & $\frac{1}2:23$; $\frac{1}6:7$ \\
11021 & $\left(\begin{smallmatrix} 1&2&3\\2&3&4\\3&4&5 \end{smallmatrix}\right)$ & TTJJ & $\frac{1}4:7$ \\
11106 & $\left(\begin{smallmatrix} 1&2&3\\2&3&4\\3&4&5 \end{smallmatrix}\right)$ & TTJJ & $\frac{1}2:15$; $\frac{1}5:6$ \\
11125 & $\left(\begin{smallmatrix} 1&2&3\\2&3&4\\2&3&4 \end{smallmatrix}\right)$ & TTJJ & $\frac{1}2:14$; $\frac{1}3:7$; $\frac{1}4:6$ \\
11455 & $\left(\begin{smallmatrix} 1&2&2\\2&3&3\\2&3&3 \end{smallmatrix}\right)$ & TTJJ & $\frac{1}2:11$; $\frac{1}3:6$ \\
16228 & $\left(\begin{smallmatrix} 1&2&3\\1&2&3\\2&3&4 \end{smallmatrix}\right)$ & TTJJ & $\frac{1}2:9$; $\frac{1}4:5$ \\
16339 & $\left(\begin{smallmatrix} 1&2&2\\1&2&2\\2&3&3 \end{smallmatrix}\right)$ & TTJJ & $\frac{1}2:8$; $\frac{1}3:5$ \\
20652 & $\left(\begin{smallmatrix} 1&2&2\\1&2&2\\1&2&2 \end{smallmatrix}\right)$ & TTJ & $\frac{1}2:6$ \\
24078 & $\left(\begin{smallmatrix} 1&1&2\\1&1&2\\2&2&3 \end{smallmatrix}\right)$ & TTJ & $\frac{1}3:4$ \\
\end{longtable}

For example, for Hilbert series number 4839,
\[
P_{4839}(t) = \frac{1 - t^{11} - 2t^{12} - 2t^{13} - 2t^{14} - t^{15} - t^{16} + \cdots - t^{40}
}{\prod_{a \in [1,1,4,5,6,7,8,9]} (1-t^a)},
\]
\cite{BKRbigtable} describes 4 deformation families of Fano 3-folds
\[
X\subset\P(1,1,4,5,6,7,8,9).
\]
A general such $X$ has Type~I projections from both $\frac{1}5(1,1,4)$
and $\frac{1}9(1,1,8)$.
(It is enough to consider just one of these centres of projection, but
\cite{TJ} calculates both, drawing the same conclusion twice.)

We construct a $\PxP$ model for $P_{4839}$.
Consider $\P=\P^7(1,1,4,5,6,7,8,9)$ with coordinates $x,y,z,t$, $u,v,w,s$.
The $2\times 2$ minors of the matrix
\[
\begin{pmatrix}
     z &  u &  v \\
    t &  v + x^7-y^7 & w + z^2 + x^8 \\
    u + x^6+y^6 & w & s
\end{pmatrix}
\quad\text{of weights}\quad
\begin{pmatrix}
 4 & 6 & 7 \\ 
 5 & 7 & 8 \\
 6 & 8 & 9
\end{pmatrix}
\]
define quasismooth $X\subset\P$ with quotient singularities
$\frac{1}2(1,1,1)$, $\frac{1}5(1,1,4)$ and $\frac{1}9(1,1,8)$.

Eliminating either the variable $t$ of degree~5 or $s$ of degree~9
computes the two possible Type~I projections,
with image a nodal codimension~3 Fano 3-fold $Y$ containing
$D=\P(1,1,4)$ or $D=\P(1,1,8)$ with 20 or 13 nodes lying on $D$ respectively.
(Both $t$ and $s$ appear only once in the matrix, so eliminating them
simply involves omitting that entry and mounting the rest of the matrix
in a skew matrix, as usual.)

\section{Cases with no numerical Type I projection}
\label{s!noproj}

The five Hilbert series 360, 577, 648, 878 and 1766 
do not admit a Type~I projection, and so the analysis of \cite{TJ}
does not apply. Nevertheless each is realised by a variety in $\PxP$
format exist, although only two of these are Fano 3-folds.

In the two cases 360 and 648 the general $\PxP$ model
is not quasismooth and has a non-terminal singularity,
so there is no $\PxP$ Fano model.
(Each of these admit Type~II$_1$ projections,
so are instead subject to the analysis of \cite{Pap}; this is 
carried out by Taylor~\cite{RT}.)
In the case 577, the $\PxP$ model is quasismooth, but it
has a $\frac{1}4(1,1,1)$ quotient singularity and so is not
a terminal Fano 3-fold and again there is no $\PxP$ Fano model.

However, there is a quasismooth Fano 3-fold $X\subset\P(1,3^2,4^2,5^2,6)$
in $\PxP$ format with weights
\[
\begin{pmatrix} 3&4&5\\4&5&6\\5&6&7 \end{pmatrix}
\]
realising Hilbert series 878.
It has $4\times\frac{1}3(1,1,2)$, $2\times\frac{1}4(1,1,3)$
quotient singularities.
There is also a quasismooth Fano 3-fold $X\subset\P(1,2,3^3,4^2,5)$
in $\PxP$ format with weights
\[
\begin{pmatrix} 2&3&4\\3&4&5\\4&5&6 \end{pmatrix}
\]
realising Hilbert series 1799.
It has $2\times\frac{1}2(1,1,1)$, $5\times\frac{1}3(1,1,2)$
quotient singularities.
Each of these two admit only Type~II$_2$ projections, and
an analysis by Gorenstein projection has not yet been attempted.
Presumably such an analysis can in principle work, once we have
much better understanding of Type~II unprojection, but until then our models
are the only Fano 3-folds known to
realise these two Hilbert series.

\section{Enumerating $\PxP$ formats}
\label{s!formats}

\subsection{Enumerating $\PxP$ formats and cases that fail}
\label{s!enum}

The Hilbert series $P_{X}(t) = \sum_{m\in\N} h^0(-mK_X) t^m$
of such Gorenstein rings $R(X)$ satisfy the
orbifold integral plurigenus formula \cite[Theorem~1.3]{BRZ}
\begin{equation}\label{e!icecream}
P_X(t) = \Pini(t) + \sum_{Q\in\cB} \Porb(Q)(t),
\end{equation}
where $\Pini$ is a function only of the genus $g$ of $X$, where $g+2=h^0(-K_X)$,
and $\Porb$ is a function of a quotient singularity $Q = \frac{1}r(1,a,-a)$,
the collection of which form the basket $\cB$ of $X$ (see \cite[\S9]{CPR}).
When $X\subset \wP$ is quasismooth, and so is an orbifold,
the basket $\cB$ is exactly the collection of quotient singularities of $X$.
Thus the numerical data $g, \cB$ gives the basis for a systematic search
of Hilbert series with given properties, which we develop further here.

We may enumerate all $\PxP$ formats $V(a,b)$ and then list
all genus--basket pairs $g,\cB$ whose corresponding
series \eqref{e!icecream} has matching numerator.
This algorithm is explained in \cite[\S4]{BKZ}.
It works systematically through increasing $k\in\N$,
where $k=3(\sum a_i + \sum b_i)$, the sum of the weights
of the ambient space of the image of $\Phi$ in \eqref{eq!phi}.

The enumeration does not have a termination condition,
even though there can only be finitely many solutions for Fano 3-folds,
so this does not directly give a classification.
Nevertheless, we search for $\P^2\times\P^2$ formats for each $k=1,\dots,31$
to start the investigation.
This reveals 53 cases whose numerical data (basket and genus) match
those of a Fano 3-fold.
The number \#\ of cases found per value of $k$ is:
\[
\begin{array}{c|ccccccccccccccccccccc}
k & 4 & 5 & 6 & 7 & 8 & 9 & 10 & 11 & 12 & 13 & 14 & 15 & 16 & 17 & 18 & 19 & 20 & 21 & 22 & 23 & \hbox{24--31} \\
\hline
\# & 1 & 3 & 2 & 3 & 3 & 5 & 4 & 4 & 4 & 3 & 5 & 5 & 1 & 2 & 3 & 1 & 1 & 0 & 2 & 1 & 0
\end{array}
\]
This hints that we may have found all Fano Hilbert series that
match some $\PxP$ format, since the algorithm stops producing
results after $k=23$. Of course that is not a proof that there
are no other cases, and we do not claim that;
the results here only use the outcome of this search as their
starting point, so how that outcome arises is not relevant.

\subsection{Weighted ${\mathrm{GL}}(3,\C)\times {\mathrm{GL}}(3,\C)$ varieties according to Szendr{\H o}i}
\label{s!sz05}
The elementary considerations we deploy for the key varieties $V(a,b)$
are part of a more general approach to weighted homogeneous spaces
by Grojnowski and Corti--Reid \cite{CR}, with
other cases developed by Qureshi and Szendr{\H o}i \cite{QS11,QS}.
The particular case of $\PxP$ was worked out detail by Szendr{\H o}i
\cite{sz05}, which we sketch here.

In the treatment of \cite{sz05}, $G = \GL(3,\C)\times \GL(3,\C)$
has weight lattice $M=\Hom(T,\C^*)\cong\Z^6$, for the maximal torus $T\subset G$.
The construction of a weighted $\PxP$, denoted $w\Sigma(\mu,u)$,
is determined by the choice
of a coweight vector $\mu\in \Hom(M,\Z)$, in coordinates say
$\mu = (a_1,a_2,a_3,b_1,b_2,b_3)\in \Hom(M,\Z)$, and an integer $u\in\Z$.
These data are subject to the positivity conditions that all $a_i+b_j+u>0$.
The construction of $w\Sigma(\mu,u)$ is described in \cite[\S2.2]{QS11}.
It embeds in wps
\begin{equation}\label{e!wSigma}
w\Sigma(\mu,u) \hookrightarrow w\P^8(a_1+b_1+u,\dots,a_3+b_3+u),
\end{equation}
with image defined by $2\times 2$ minors
\begin{equation}\label{e!image}
w\Si=\left\{\left.\bigwedge^2\left(\begin{array}{ccc}
x_1 & x_2 & x_3 \\
x_4 & x_5 & x_6 \\
x_7 & x_8 & x_9 \\
\end{array}\right)\right.= 0 \right\} \subset w\P^8
\end{equation}
with respect to the weights
\[
\deg
\left(\begin{array}{ccc}
x_1 & x_2 & x_3 \\
x_4 & x_5 & x_6 \\
x_7 & x_8 & x_9 \\
\end{array}\right)
=
\left(\begin{array}{ccc}
a_1+b_1+u & a_1+b_2+u & a_1+b_3+u \\
a_2+b_1+u & a_2+b_2+u & a_2+b_3+u \\
a_3+b_1+u & a_3+b_2+u & a_3+b_3+u\\
\end{array}\right).
\]
The following theorem then follows from the general Hilbert series
formula of \cite[Theorem 3.1]{QS11}.

\begin{thm}[Szendr{\H o}i \cite{sz05}]
The Hilbert series of $w\Sigma(\mu,u)$ in the embedding \eqref{e!wSigma} is
\[
P(t) =
\dfrac{\Pnum(t)}{\prod_{i,j}\left(1-t^{a_i+b_j+u}\right)},
\]
where the Hilbert numerator $\Pnum(t)$ is
\[
1-\left(\sum_{i,j}t^{-a_i-b_j}\right)t^{2u+s}+\left(4+\sum_{i\ne j}t^{-a_i+a_j}+\sum_{i\ne j}t^{-b_i+b_j}\right)t^{3u+s} -\left(\sum_{i,j}t^{a_i+b_j}\right)t^{4u+s}+t^{6u+2s},
\]
with $s=a_1+a_2+a_3+b_1+b_2+b_3$.
\end{thm}
This numerator exposes the $9\times 16$ resolution.
The $2\times 2$ minors in \eqref{e!image}
are visible in the first parentheses; for example $t^{-a_1-b_1}t^{2u+s} = t^{(a_2+b_2+u) + (a_3+b_3+u)}$ carries the degree of $x_5x_9=x_6x_8$.
First syzygies appear in the second parentheses; for example, the syzygy
\[
\det
\left(\begin{array}{ccc}
x_4 & x_5 & x_6 \\
x_4 & x_5 & x_6 \\
x_7 & x_8 & x_9 \\
\end{array}\right)
\equiv 0
\]
has degree
$\deg(x_4x_5x_9) =  (a_2+b_1+u) + (a_2+b_2+u) + (a_3+b_3+u) = (a_2-a_1) + 3u + s$.
The additional parameter $u\in\Z$ in this treatment is absorbed into
the $a_i$ in our naive treatment of \S\ref{s!V}, so the key varieties
we enumerate are the same.

\subsection*{Acknowledgments}

It is our pleasure to thank Bal{\'a}zs Szendr{\H o}i for sharing his unpublished
analysis \cite{sz05} of weighted $\GL(3,\C)\times\GL(3,\C)$ varieties that provided
our initial motivation and the tools of \S\ref{s!sz05}, and for several helpful 
and informative conversations.
AMK was supported by EPRSC Fellowship EP/N022513/1.
MIQ was supported by a Lahore University of Management Sciences (LUMS)
faculty startup research grant (STG-MTH-1305)
and an LMS research in pairs grant for a visit to the UK.

\bibliographystyle{alpha}

\begin{thebibliography}{KMMT00}

\bibitem[ABR02]{ABR02}
Selma Alt{\i}nok, Gavin Brown, and Miles Reid.
\newblock Fano 3-folds, {$K3$} surfaces and graded rings.
\newblock In {\em Topology and geometry: commemorating {SISTAG}}, volume 314 of
  {\em Contemp. Math.}, pages 25--53. Amer. Math. Soc., Providence, RI, 2002.

\bibitem[BF17]{BF}
Gavin Brown and Enrico Fatighenti.
\newblock Hodge numbers and deformations of {F}ano 3-folds.
\newblock Preprint available online as arXiv:1707.00653, 2017.

\bibitem[BK]{grdb}
Gavin Brown and Alexander~M. Kasprzyk.
\newblock The graded ring database.
\newblock Online.
\newblock Access via
  \href{http://www.grdb.co.uk/}{\texttt{http://www.grdb.co.uk/}}.

\bibitem[BKR12a]{TJ}
Gavin Brown, Michael Kerber, and Miles Reid.
\newblock Fano 3-folds in codimension 4, {T}om and {J}erry. {P}art {I}.
\newblock {\em Compos. Math.}, 148(4):1171--1194, 2012.

\bibitem[BKR12b]{BKRbigtable}
Gavin Brown, Michael Kerber, and Miles Reid.
\newblock {T}om and {J}erry: {B}ig {T}able, 2012.
\newblock Linked at
  \href{http://grdb.co.uk/Downloads/}{\texttt{http://grdb.co.uk/Downloads/}}.

\bibitem[BKZ14]{BKZ}
Gavin Brown, Alexander~M. Kasprzyk, and Lei Zhu.
\newblock Canonical threefolds in {G}orenstein formats.
\newblock Preprint available online as arXiv:1409.4644, 2014.

\bibitem[BR13]{dip1}
Gavin Brown and Miles Reid.
\newblock Diptych varieties, {I}.
\newblock {\em Proc. Lond. Math. Soc. (3)}, 107(6):1353--1394, 2013.

\bibitem[Bro06]{k3degen}
Gavin Brown.
\newblock Graded rings and special {$K3$} surfaces.
\newblock In {\em Discovering mathematics with {M}agma}, volume~19 of {\em
  Algorithms Comput. Math.}, pages 137--159. Springer, Berlin, 2006.

\bibitem[Bro07]{k3db}
Gavin Brown.
\newblock A database of polarized {$K3$} surfaces.
\newblock {\em Experiment. Math.}, 16(1):7--20, 2007.

\bibitem[BRZ13]{BRZ}
A.~Buckley, M.~Reid, and S.~Zhou.
\newblock Ice cream and orbifold {R}iemann-{R}och.
\newblock {\em Izv. Ross. Akad. Nauk Ser. Mat.}, 77(3):29--54, 2013.

\bibitem[CCC11]{ccc}
Jheng-Jie Chen, Jungkai~A. Chen, and Meng Chen.
\newblock On quasismooth weighted complete intersections.
\newblock {\em J. Algebraic Geom.}, 20(2):239--262, 2011.

\bibitem[CPR00]{CPR}
Alessio Corti, Aleksandr Pukhlikov, and Miles Reid.
\newblock Fano {$3$}-fold hypersurfaces.
\newblock In {\em Explicit birational geometry of 3-folds}, volume 281 of {\em
  London Math. Soc. Lecture Note Ser.}, pages 175--258. Cambridge Univ. Press,
  Cambridge, 2000.

\bibitem[CR02]{CR}
Alessio Corti and Miles Reid.
\newblock Weighted {G}rassmannians.
\newblock In {\em Algebraic geometry}, pages 141--163. de Gruyter, Berlin,
  2002.

\bibitem[DNFF15]{DFF}
C.~Di~Natale, E.~Fatighenti, and D.~Fiorenza.
\newblock Hodge theory and deformations of affine cones of subcanonical
  projective varieties.
\newblock Preprint available online as arXiv:1512.00835, 2015.

\bibitem[GS]{M2}
Daniel~R. Grayson and Michael~E. Stillman.
\newblock Macaulay2, a software system for research in algebraic geometry.
\newblock Available at \url{http://www.math.uiuc.edu/Macaulay2/}.

\bibitem[GW78]{GW}
Shiro Goto and Keiichi Watanabe.
\newblock On graded rings. {I}.
\newblock {\em J. Math. Soc. Japan}, 30(2):179--213, 1978.

\bibitem[IF00]{Fletcher}
A.~R. Iano-Fletcher.
\newblock Working with weighted complete intersections.
\newblock In {\em Explicit birational geometry of 3-folds}, volume 281 of {\em
  London Math. Soc. Lecture Note Ser.}, pages 101--173. Cambridge Univ. Press,
  Cambridge, 2000.

\bibitem[Ilt12]{ilten}
Nathan~Owen Ilten.
\newblock Versal deformations and local {H}ilbert schemes.
\newblock {\em J Softw. Algebra Geom.}, (4):12--16, 2012.

\bibitem[Kaw92]{K92}
Yujiro Kawamata.
\newblock Boundedness of {$\bold Q$}-{F}ano threefolds.
\newblock In {\em Proceedings of the {I}nternational {C}onference on {A}lgebra,
  {P}art 3 ({N}ovosibirsk, 1989)}, volume 131 of {\em Contemp. Math.}, pages
  439--445. Amer. Math. Soc., Providence, RI, 1992.

\bibitem[KM83]{KM}
A.~Kustin and M.~Miller.
\newblock Constructing big {G}orenstein ideals from small one.
\newblock {\em J. Alg}, 85:303--322, 1983.

\bibitem[KMMT00]{kmm}
J{\'a}nos Koll{\'a}r, Yoichi Miyaoka, Shigefumi Mori, and Hiromichi Takagi.
\newblock Boundedness of canonical {$\bold Q$}-{F}ano 3-folds.
\newblock {\em Proc. Japan Acad. Ser. A Math. Sci.}, 76(5):73--77, 2000.

\bibitem[Pap04]{Pwith}
Stavros~Argyrios Papadakis.
\newblock Kustin-{M}iller unprojection with complexes.
\newblock {\em J. Algebraic Geom.}, 13(2):249--268, 2004.

\bibitem[Pap08]{Pap}
Stavros~Argyrios Papadakis.
\newblock The equations of type {$\rm II\sb 1$} unprojection.
\newblock {\em J. Pure Appl. Algebra}, 212(10):2194--2208, 2008.

\bibitem[PR04]{PR}
Stavros~Argyrios Papadakis and Miles Reid.
\newblock Kustin-{M}iller unprojection without complexes.
\newblock {\em J. Algebraic Geom.}, 13(3):563--577, 2004.

\bibitem[QS11]{QS11}
Muhammad~Imran Qureshi and Bal{\'a}zs Szendr{\H o}i.
\newblock Constructing projective varieties in weighted flag varieties.
\newblock {\em Bull. Lond. Math. Soc.}, 43(4):786--798, 2011.

\bibitem[QS12]{QS}
Muhammad~Imran Qureshi and Bal{\'a}zs Szendr{\H o}i.
\newblock Calabi-{Y}au threefolds in weighted flag varieties.
\newblock {\em Adv. High Energy Phys.}, pages Art. ID 547317, 14, 2012.

\bibitem[Qur15]{qureshi}
Muhammad~Imran Qureshi.
\newblock Constructing projective varieties in weighted flag varieties {II}.
\newblock {\em Math. Proc. Cambridge Philos. Soc.}, 158(2):193--209, 2015.

\bibitem[Qur17]{QJSC}
Muhammad~Imran Qureshi.
\newblock Computing isolated orbifolds in weighted flag varieties.
\newblock {\em Journal of Symbolic Computation}, 79, Part 2:457--474, 2017.

\bibitem[Rei11]{fun}
Miles Reid.
\newblock Fun in codimension~$4$.
\newblock Preprint available online via the author's webpage, 2011.

\bibitem[Rei15]{reid4}
Miles Reid.
\newblock Gorenstein in codimension 4: the general structure theory.
\newblock In {\em Algebraic geometry in east {A}sia---{T}aipei 2011}, volume~65
  of {\em Adv. Stud. Pure Math.}, pages 201--227. Math. Soc. Japan, Tokyo,
  2015.

\bibitem[Sze05]{sz05}
Bal{\'a}zs Szendr{\H o}i.
\newblock On weighted homogeneous varieties.
\newblock Unpublished manuscript, 2005.

\bibitem[Tak02]{takagi02}
Hiromichi Takagi.
\newblock On classification of {$\Bbb Q$}-{F}ano 3-folds of {G}orenstein index
  2. {I}, {II}.
\newblock {\em Nagoya Math. J.}, 167:117--155, 157--216, 2002.

\bibitem[Tay]{RT}
Rosemary Taylor.
\newblock Fano 3-folds and {T}ype {II} {U}nprojection, {P}h{D} thesis,
  {U}niversity of {W}arwick.
\newblock In preparation.

\end{thebibliography}

\noindent {\sc Gavin Brown, Mathematics Institute, University of Warwick, UK}

\noindent {\tt G.Brown@warwick.ac.uk}

\vspace{0.1in}

\noindent {\sc Alexander M. Kasprzyk, School of Mathematical Sciences, University of Nottingham, UK}

\noindent {\tt a.m.kasprzyk@nottingham.ac.uk}

\vspace{0.1in}

\noindent {\sc Muhammad Imran Qureshi, Department of Mathematics, SBASSE, LUMS, Pakistan}

\noindent {\tt imran.qureshi@lums.edu.pk}
\end{document}